\documentclass{amsart}

\usepackage{a4}

\newtheorem{theorem}{Theorem}[section]
\newtheorem{lemma}[theorem]{Lemma}

\newtheorem{corollary}[theorem]{Corollary}

\newtheorem{example}[theorem]{Example}
\def\N{\mathbf N}
\def\RP{{\bf RP}}
\def\CP{{\bf CP}}
\def\cf{cf.}
\makeatletter
\newenvironment{proofnoqed}[1][\proofname]{\par \normalfont
  \topsep6\p@\@plus6\p@ \trivlist
  \item[\hskip\labelsep\it
    #1\@addpunct{.}]\ignorespaces
}{\endtrivlist}
\def\qedhere#1{\hbox to0pt{\hskip#1$\qed$\hss}}
\makeatother

\newtheorem{fact}[theorem]{Fact}
\DeclareMathOperator{\Sq}{Sq}
\DeclareMathOperator{\SQ}{SQ}
\def\A{\mathcal A}
\def\F{\mathbf F}

\title{Steenrod squares on conjugation spaces}
\author{Matthias Franz and Volker Puppe}
\email{matthias.franz@ujf-grenoble.fr}
\email{volker.puppe@uni-konstanz.de}
\address{Fachbereich Mathematik, 
  Universit\"at Konstanz, 78457 Konstanz, Germany}

\begin{document}

\begin{abstract}
  We prove that the coefficients of the so-called conjugation equation
  for conjugation spaces in the sense of Hausmann--Holm--Puppe are completely
  determined by Steenrod squares. This generalises a result of V.\,A.~Krasnov
  for certain complex algebraic varieties. It also leads to a generalisation of
  a formula given by Borel and Haefliger, thereby largely answering
  an old question of theirs in the affirmative.
\end{abstract}

\maketitle

\section{Statement of the results}
\noindent
Let $X$ be a topological space with an involution~$\tau$.
We look at $X$ as a space with an action of the group~$G=\{1,\tau\}$.
We take cohomology with coefficients in~$\F_2$
and consider the restriction map~$r\colon H^*(X)\to H^*(X^\tau)$,
its equivariant 
counterpart $r_G\colon H_G^*(X)\to H_G^*(X^\tau)=H^*(X^\tau)\otimes\F_2[u]$
and the canonical projection~$p\colon H_G^*(X)\to H^*(X)$.

According to Hausmann--Holm--Puppe~\cite{HausmannHolmPuppe:05},
$(X,\tau)$ is called a \emph{conjugation space}
if there exists
a section~$\sigma\colon H^*(X)\to H_G^*(X)$
to~$p$ 
and a degree-halving isomorphism~$\kappa\colon H^*(X)\to H^*(X^\tau)$
with the following property: for every~$x\in H^n(X)$, $n\in\N$,
there exists elements~$y_1$,~\ldots,~$y_n\in H^*(X^\tau)$
such that the so-called \emph{conjugation equation}
\begin{equation}
  r_G(\sigma(x))=\kappa(x)u^n+y_1 u^{n-1}+\cdots+y_{n-1}u+y_n
\end{equation}
holds.
A priori, $\sigma$ and $\kappa$ are only assumed to be additive,
but the conjugation equation implies that they are in fact multiplicative
and unique.

The coefficients~$y_i$ in the conjugation equation are 
completely determined by~$x$. We prove that they are the Steenrod squares
of~$\kappa(x)$.

\begin{theorem}\label{r-sigma}
  For every~$x\in H^{2n}(X)$, $n\in\N$, one has
  $$
    r_G(\sigma(x))=\sum_{i=0}^n\Sq^i(\kappa(x))\,u^{n-i}.
  $$
\end{theorem}

\begin{corollary}\label{square-kappa}
  For every~$x\in H^*(X)$, one has
  $$
    r(x)=\kappa(x)^2.
  $$
\end{corollary}

We also show that the isomorphism~$\kappa$
commutes with total Steenrod squares.

\begin{theorem}\label{kappa-Sq}
  For every $x\in H^*(X)$ one has
  $$
    \kappa(\Sq(x))=\Sq(\kappa(x)).
  $$
\end{theorem}

Note that the odd Steenrod squares of~$x$ vanish
since $H^*(X)$ is concentrated in even degrees.
Hence, the above identity is equivalent to
\begin{equation}
  \kappa(\Sq^{2k}(x))=\Sq^k(\kappa(x))
  \quad\hbox{for all~$k\in\N$.}
\end{equation}

\section{Proofs}
\noindent
We denote the Steenrod algebra for the prime~$2$ by~$\A$.

\begin{lemma}\label{operations-are-squares}
  For every $n$ there exist elements~$a_0$,~\ldots,~$a_n$,~$b\in\A$
  such that for every conjugation space~$X$ and every~$x\in H^{2n}(X)$ one has
  $$
    r_G(\sigma(x))=\sum_{i=0}^n a_i(\kappa(x))\,u^{n-i}
    \qquad\text{and}\qquad
    \kappa(\Sq(x))=b(\kappa(x)).
  $$
  Moreover, $a_0=1$ and $a_1=\Sq^1$.
\end{lemma}

\begin{proofnoqed}
  Since $\kappa$ is bijective, one can define functions~$a_i$ and~$b$
  for each~$X$
  such that the above identities hold. We show that they are
  (or, more precisely, come from) Steenrod squares, using that
  the restriction map~$r_G$ commutes with all Steenrod squares.
  We write $\kappa(x)=z$.

  We start by proving the claim about the~$a_i$ by induction,
  beginning at~$a_0(z)=z$.
  If $i>0$ is even, we apply $\Sq^{2k}$, where $k\le i/2$ will be chosen later.
  We can write
  \begin{equation}\label{Sq-2k-sigma-x}
    \Sq^{2k}(\sigma(x))=\sum_{l=-n}^k\sigma(\tilde x_l)\,u^{2(k-l)}
  \end{equation}
  for some~$\tilde x_l\in H^{2(n+l)}(X)$.
  Write $z_l=\kappa(\tilde x_l)$.
  The restriction~$r_G(\sigma(\tilde x_l)\,u^{2(k-l)})$
  has leading term~$z_l u^{n+2k-l}$,
  while the leading power of~$u$ in~$\Sq^{2k}(r_G(\sigma(x)))$
  is at most~$u^{n+2k}$. Hence, the summation in~\eqref{Sq-2k-sigma-x}
  is in fact only over~$0\le l\le k$.

  Comparing coefficients of~$u^{n+2k-l}$ yields 
  \begin{equation}\label{formula-z-l}
    z_l=\sum_{j=0}^l\binom{n-l+j}{2k-j}\Sq^j(a_{l-j}(z))
        +\sum_{j=1}^l\Sq^j(z_{l-j}).
  \end{equation}
  This shows $z_l=b_l(z)$ for some~$b_l\in\A$.
  Comparing coefficients of~$u^{n+2k-i}$ gives
  \begin{align*}
    \sum_{l=0}^k a_{i-l}(z_l)
    &= \binom{n}{2k}\,a_i(z)+\sum_{l=1}^k a_{i-l}(b_l(z)) \\
    &= \binom{n-i}{2k}\,a_i(z)
        +\sum_{j=1}^{2k}\binom{n-i+j}{2k-j}\Sq^j(a_{i-j}(z)).
  \end{align*}

  Now suppose that $k\le i/2$ is such that
  $$
    \binom{n}{2k}\ne\binom{n-i}{2k}.
  $$
  For instance, this is true if $2k$ is the largest power of~$2$ dividing~$i$.
  (Recall that a binomial coefficient mod~$2$ is the product of the binomial
  coefficients taken for each pair of binary digits,
  \cf~\cite[Lemma~I.2.6]{SteenrodEpstein:62}.)
  Then the above equation can be solved for~$a_i(z)$ and shows
  that $a_i(z)$ can be expressed in terms of repeated Steenrod squares of~$z$.

  For odd~$i$, a similar (but simpler) reasoning based on
  commutativity with respect to~$\Sq^1$ gives
  $a_i(z)=\Sq^1(a_{i-1}(z))$, in particular $a_1(z)=\Sq^1(z)$.

  Now that all~$a_i(z)$ are known, we apply $\Sq^{2k}$ for any~$k$.
  Using the same notation as above, we have
  $\Sq^{2k}(x)=p(\Sq^{2k}(\sigma(x)))=\tilde x_k$.
  Comparing coefficients as before gives a formula
  for~$b_l(z)$ similar to equation~\eqref{formula-z-l},
  but where the summation index~$j$ starts at~$l-n$ if~$l>n$.
  Still, the equations can be recursively solved for~$z_l$.
  Hence,
  $$
    \kappa(\Sq(x))=\kappa(\tilde x_0)+\cdots+\kappa(\tilde x_n)
      =b_0(z)+\cdots+b_n(z)=b(z).
    \qedhere{32pt}
  $$
\end{proofnoqed}

In principle, the preceding proof could be used to determine
the coefficients of the conjugation equation completely
(as well as those of~$\Sq^k(\sigma(x))$ for any~$k$). We will take
a less tedious approach which relies on the fact that suitable
products of infinite-dimensional real projective space can
``detect'' Steenrod squares, \cf~\cite[Corollary~I.3.3]{SteenrodEpstein:62}.

\begin{fact}\label{Sq-injective}
  The restricted evaluation map
  $$
    \A_{\le n}\to H^*((\RP^\infty)^n),
    \quad
    a\mapsto a(u\times\cdots\times u)
  $$
  is injective for any~$n\in\N$.
\end{fact}


\begin{proof}[Proof of Theorem~\ref{r-sigma}]
  Write $\SQ$ for the homogenised total Steenrod square,
  $$
    \SQ(\kappa(x))=\sum_{i=0}^n\Sq^i(\kappa(x))\,u^{n-i}.
  $$
  We want to show $r_G(\sigma(x))=\SQ(\kappa(x))$
  for all cohomology classes of all conjugation spaces
  By Lemma~\ref{operations-are-squares} and
  Fact~\ref{Sq-injective}, it suffices
  to do so for the conjugation space~$X=(\CP^\infty)^n$
  with complex conjugation and $x$, the $n$-fold cross product
  of the generator~$v\in H^2(\CP^\infty)$ because $X^\tau=(\RP^\infty)^n$
  and $\kappa(x)=u\times\cdots\times u$ in this case.
  For~$n=1$ the identity is true since we already know $a_1$.
  The general case reduces to the case~$n=1$ because of the multiplicativity
  of all maps involved:
  \begin{align*}
    r_G(\sigma(x)
    &= r_G(\sigma(v\times\cdots\times v))
     = r_G(\sigma(v))\times\cdots\times r_G(\sigma(v)) \\
    &= \SQ(\kappa(v))\times\cdots\times \SQ(\kappa(v))
     = \SQ(\kappa(v\times\cdots\times v))
     = \SQ(\kappa(x)).
  \end{align*}
  (Note that $\kappa$,~$\sigma$ and~$r_G$ commute with products.)
\end{proof}

\begin{proofnoqed}[Proof of Corollary~\ref{square-kappa}]
  We have for~$x\in H^{2n}(X)$
  $$
    r(x)=r(p(\sigma(x)))=p(r_G(\sigma(x)))=p(\Sq^n(\kappa(x)))=\kappa(x)^2.
    \qedhere{37.5pt}
  $$
\end{proofnoqed}

\begin{proofnoqed}[Proof of Theorem~\ref{kappa-Sq}]
  As in the proof of Theorem~\ref{r-sigma}, it suffices to show
  the claimed identity
  for~$X=(\CP^\infty)^n$ and $x=v\times\cdots\times v$.
  Again, the general case can be reduced to~$n=1$, where we find
  $$
    \kappa(\Sq(v))=\kappa(v+v^2)=\kappa(v)+\kappa(v)^2=\Sq(\kappa(v)).
    \qedhere{56.5pt}
  $$
\end{proofnoqed}

\section{Remarks}
\noindent
  Let $X$ be a non-singular complex projective variety
  defined over the reals
  such that its real locus~$X^\tau$ is non-empty.
  All algebraic cycles in~$X$ are understood to be defined
  over the reals.
  Borel--Haefliger have shown that
  if $H_*(X)$~and~$H_*(X^\tau)$ are generated by algebraic cycles,
  then the restriction~$\lambda$ of cycles in~$X$ to their real locus
  induces a degree-halving isomorphism~$H_*(X)\to H_*(X^\tau)$
  respecting intersection products \cite[\S 5.15]{BorelHaefliger:61}.
  They also show that if $H^*(X^\tau)$ is generated by algebraic cycles
  and $x\in H^*(X)$ is Poincar\'e dual to a linear combination of non-singular
  subvarieties, then the identity in Theorem~\ref{kappa-Sq} holds, and they ask
  whether it holds more generally \cite[\S 5.17]{BorelHaefliger:61}.

  V.\,A.~Krasnov has proved that for a variety~$X$ as above,
  Theorem~\ref{r-sigma} holds for cohomology classes Poincar\'e dual
  to algebraic cycles, where $\kappa$
  is the Poincar\'e transpose of~$\lambda$
  and $\sigma$~the canonical section \cite[Theorem~4.2]{Krasnov:94}.
  This implies that if $H_*(X)$ is generated by algebraic cycles,
  then so is $H_*(X^\tau)$ \cite[Theorem~0.1]{Krasnov:03}.
  Moreover, $X$ is a  conjugation spaces in the sense
  of~\cite{HausmannHolmPuppe:05}.

  In a topological framework van~Hamel has recently shown
  that certain topological manifolds with involutions are conjugation spaces
  \cite[Theorem]{vanHamel:??}. The necessary assumptions
  are formulated in terms of topological cycles.

  The following simple example shows that in general the existence
  of a degree-halving
  multiplicative isomorphism~$\kappa\colon H^*(X)\to H^*(X^\tau)$
  by itself does not imply that $(X,\tau)$ is a conjugation space.

\begin{example}\rm
  Let $X=S^2\times S^4$ be equipped with the componentwise involution~$\tau$
  which is the identity for~$S^2$ and for~$S^4$ has fixed point set~$S^1$.
  So $X^\tau=S^2\times S^1$. Clearly there is a degree-halving multiplicative
  isomorphism~$\kappa\colon H^*(X)\to H^*(X^\tau)$.
  It is easy to check there is also a multiplicative
  section~$\sigma\colon H^*(X)\to H_G^*(X)$.
  But $(X,\tau)$ is not a conjugation space:
  the restriction map
  $$
    r_G\colon H_G^*(S^2\times S^4)\cong H^*(S^2\times S^4)\otimes \F_2[u]
    \to H^*(S^2\times S^1)\otimes \F_2[u]
  $$
  is given by
  $s_2\otimes1\mapsto s_2\otimes1$ and $s_4\otimes1\mapsto s_1\otimes u^3$,
  where $s_n\in H^n(S^n)$ denotes the generator. Hence the conjugation equation
  does not hold. Of course, $S^2\times S^4$ with the different componentwise
  involution~$\tilde\tau$ which has $S^1\subset S^2$ and $S^2\subset S^4$
  as fixed point sets (and hence $X^{\tilde\tau}=S^1\times S^2\cong X^\tau$)
  is a conjugation space.
\end{example}

\bigbreak

\end{document}